\documentclass[letterpaper, 10 pt, conference]{ieeeconf}  % Comment this line out if you need a4paper

\IEEEoverridecommandlockouts                              % This command is only needed if
\usepackage{amsmath}
\usepackage{amssymb}
\usepackage{bbold}
\usepackage[pdftex,final]{graphicx}

\newtheorem{prop}{Proposition}
\newtheorem{rem}{Remark}
\newtheorem{dfn}{Definition}

\newtheorem{corollary}{Corollary}

\numberwithin{equation}{section}

\def\theequation{\arabic{section}.\arabic{equation}}

\IEEEoverridecommandlockouts                             
\title{\LARGE \bf{Sufficient Lie Algebraic Conditions for Sampled-Data Feedback Stabilization}}

\author{J.~Tsinias$^{1}$ and D.~Theodosis$^{2}$% <-this % stops a space
\thanks{$^{1}$John Tsinias is with the Faculty of Mathematics, Department of Mathematical and Physical Sciences,
        National Technical University of Athens, Zografou Campus 15780, Athens, Greece
        {\tt\small jtsin@central.ntua.gr}}%
\thanks{$^{2}$Dionysis Theodosis is with the Department of Mathematical and Physical Sciences,
        National Technical University of Athens, Zografou Campus 15780, Athens, Greece
        {\tt\small dtheodp@central.ntua.gr}}%
}

\begin{document}
\maketitle
\thispagestyle{empty}
\pagestyle{empty}

\begin{abstract}
For nonlinear affine in the control systems, a Lie algebraic sufficient condition for sampled-data feedback semi-global stabilization is established. We use this result, in order to derive sufficient conditions for sampled-data feedback stabilization for a couple of three-dimensional systems. 
\end{abstract}

\section{Introduction}
\label{S:1}

Significant results towards stabilization of  nonlinear systems by means of sampled-data feedback control have appeared in the literature (see for instance \cite{art:1}, \cite{art:2}, \cite{art:4}-\cite{art:7}, \cite{art7.1}-\cite{art:11b}, \cite{art:12}, \cite{art:15} and relative references therein). In the recent works \cite{art:17} and \cite{art:19}, the concept of \textit{Weak Global Asymptotic Stabilization by Sampled-Data Feedback} (SDF-WGAS) is introduced for autonomous systems:
\begin{equation}\label{1.1}
\begin{array}{c}
\dot{x}=f(x,u),\, \, (x,u)\in \mathbb{R}^{n}\times \mathbb{R}^{m},\\
f(0,0)=0 
\end{array}
\end{equation}
and Lyapunov-like sufficient characterizations of this property are examined. Particularly, in \cite[Proposition 2]{art:19}, a Lie algebraic sufficient condition for SDF-WGAS is established for the case of affine in the control single-input systems 
\begin{equation}\label{1.2}
\begin{array}{c}
 \dot{x}=f(x)+ug(x),\; \; (x,u)\in  \mathbb{R}^{n}\times\mathbb{R},\\
   f(0)=0  
   \end{array}
\end{equation} 
This condition constitutes an extension of the well-known ``Artstein-Sontag'' sufficient condition for asymptotic stabilization of systems \eqref{1.2} by means of an almost smooth feedback; (see \cite{art:3}, \cite{art:14} and \cite{art:16}).

Throughout the paper we adopt the following notations. For any pair of $C^{1} $ mappings $X: \mathbb{R}^{n} \to \mathbb{R}^{k}$, $Y: \mathbb{R}^{k} \to \mathbb{R}^{\ell } $ we denote $XY:=(DY)X$, $DY$ being the derivative of $Y$. By $[\cdot ,\cdot ]$ we denote the Lie bracket operator, namely, $[X,Y]=XY-YX$ for any pair of $C^{1} $ mappings $X,Y:\mathbb{R}^{ n} \to \mathbb{R}^{n} $. 

The precise statement of \cite[Proposition 2]{art:19} is the following. Assume that $f,g\in C^{2} $ and there exists a $C^{2} $, positive definite and proper function $V:{\mathbb R}^{n} \to {\mathbb R}^{+} $ such that the following implication holds:
\begin{equation} \label{1.3} \begin{array}{l} {(gV)(x)=0,x\ne 0} \\ {\Rightarrow \left\{\hspace{-.5em}\begin{array}{l} {either\; (fV)(x)<0,}\\\hspace{2.5em}(\mathrm{``Artstein-Sontag"\, condition}) \\ {or\; (fV)(x)=0;\; ([f,g]V)(x)\ne 0} \end{array}\right. } \end{array} \end{equation} 
Then system \eqref{1.2} is SDF-WGAS.

In  the present work, we first deal with the general case \eqref{1.1}, providing a Lyapunov characterization for a stronger version of SDF-WGAS. Particularly, Proposition 2 of our work asserts that for systems \eqref{1.1} the same Lyapunov characterization of SDF-WGAS, originally proposed in \cite{art:17} (see Assumption 1 below), implies \textit{Semi-Global Asymptotic Stabilization by means of a time-varying Sampled-Data Feedback }(SDF-SGAS), which is a stronger type of  SDF-WGAS. We exploit the result of Proposition 2 to establish in Proposition 3 a Lie algebraic sufficient condition for SDF-SGAS(WGAS) for systems \eqref{1.2}, much weaker than \eqref{1.3}. 

The result of Proposition 3 is then used, in order to study the SDF-SGAS for a couple of 3-dimensional affine in the control cases (Corollaries 1 and 2). 

The precise statements of Propositions 2 and 3 and of Corollaries 1 and 2 are given in Section II. Proofs of both corollaries are given in Section III. Detailed proofs of Propositions 2 and 3 can be found in \cite{art:arX}. For completeness, an outline of proof of Proposition 3 is also provided in Section III.

\section{Definitions and Main Results}
\label{S:2}
Consider system \eqref{1.1} and assume that $f:{\mathbb R}^{ n} \times {\mathbb R}^{ m} \to {\mathbb R}^{ n} $ is Lipschitz continuous. We denote by $x(\cdot )=x(\cdot ,s,x_{0} ,u)$ the trajectory of \eqref{1.1} with initial condition $x(s,s,x_{0} ,u)=x_{0} \in {\mathbb R}^{ n} $ corresponding to certain measurable and locally essentially bounded control $u:[s,T_{\max } )\to {\mathbb R}^{ m} $, where $T_{\max } =T_{\max } (s,x_{0} ,u)$ is the corresponding maximal existing time of the trajectory. 
\begin{dfn}
We say that system \eqref{1.1} is \textit{Weakly Globally Asymptotically Stabilizable by Sampled-Data Feedback} (SDF-WGAS), if for any constant $\tau >0$ there exist mappings $T: {\mathbb R}^{ n} \backslash \{ 0\} \to  {\mathbb R}^{+} \backslash \{ 0\} $ satisfying 
\begin{equation} \label{2.1} T(x)\leq \tau ,\, \, \, \forall x\in  {\mathbb R}^{ n} \setminus \{ 0\}  \end{equation} 
and $k(t,x;x_{0} ): {\mathbb R}^{+} \times {\mathbb R}^{ n} \times {\mathbb R}^{ n} \to  {\mathbb R}^{ m} $ such that for any fixed $(x,x_{0} )\in {\mathbb R}^{n}\times \mathbb{R}^{n} $ the map $k(\cdot ,x;x_{0} ): {\mathbb R}^{+} \to  {\mathbb R}^{ m} $ is measurable and locally essentially  bounded and such that for every $x_{0} \ne 0$ there exists a sequence of times
\begin{equation}\label{2.2} t_{1} :=0<t_{2} <t_{3} <\ldots <t_{\nu } <\ldots \; \, , \mathrm{with}\, \,  t_{\nu } \to \infty \end{equation}  
in such a way that the trajectory $x(\cdot )$ of the sampled-data closed loop system:
\begin{align}   \dot{x}=f(x,k(t,x(t_{i} );x_{0}& )),\,\, t\in [t_{i} ,\; t_{i+1} ),\; \,\, i=1,2,\ldots \nonumber \\  &x(0)=x_{0} \in  {\mathbb R}^{ n}   \label{2.3}\end{align} 
satisfies: 
\begin{equation} \label{2.4} t_{i+1} -t_{i} =T(x(t_{i} )),\; i=1,2,\ldots  \end{equation} 
and the following properties:
\begin{flalign}\text{Stability:}  &&\hspace{-2em}
\begin{array}{c}\forall \varepsilon >0\Rightarrow \exists \delta =\delta (\varepsilon )>0: |x(0)|\le \delta\\
\Rightarrow |x(t)|\leq \varepsilon ,\; \forall t\geq 0\end{array} 
&\label{2.5}
\end{flalign}
\begin{flalign}\text{Attractivity:}  && \mathop{\lim }\limits_{t\to \infty } x(t)=0,\; \forall x(0)\in  {\mathbb R}^{n}
&&\label{2.6}\end{flalign}
\noindent where $\left|x\right|$ denotes the Euclidean norm of the vector $x$.\\
\end{dfn}

Next we give the Lyapunov characterization of SDF-WGAS proposed in \cite{art:17} and \cite{art:19}, that constitutes a generalization of the concept of the \textit{control Lyapunov function} (see Definition 5.7.1 in \cite{art:13}).\\

\textbf{Assumption 1:} \textit{There exist a positive definite $C^{0} $ function $V: {\mathbb R}^{n} \to  {\mathbb R}^{+} $ and a function $a\in K$ (namely, $a(\cdot )$ is continuous, strictly increasing with $a(0)=0$) such that for every $\xi >0$ and $x_0\neq0$ there exists a constant $\varepsilon =\varepsilon(x_0)\in(0,\xi]$ and a measurable and locally essentially bounded control $u(\cdot,x_{0}) :[0,\varepsilon ]\to  {\mathbb R}^{m} $  satisfying}
\begin{subequations} \label{2.7}
\begin{equation}\label{2.7a}
 V(x(\varepsilon ,0,x_{0} ,u(\cdot,x_{0}) ))<V(x_{0} );  
\end{equation}
\begin{equation}\label{2.7b}
 V(x(s,0,x_{0} ,u(\cdot,x_{0}) ))\le a(V(x_{0} )),\; \; \forall s\in [0,\varepsilon ]  
\end{equation}
\end{subequations}

\noindent The following result was established in \cite{art:17}.\\

\begin{prop}
Under Assumption 1, system \eqref{1.1} is SDF-WGAS.\\
\end{prop}

We now present the concept of SDF-SGAS, which, as mentioned above, is a strong version of SDF-WGAS:\\

\begin{dfn} 
We say that system \eqref{1.1} is \textit{Semi-Globally Asymptotically Stabilizable by Sampled-Data Feedback}  (SDF-SGAS), if for every $R>0$ and for any given partition of times 
\begin{equation} \label{2.8} T_{1} :=0<T_{2} <T_{3} <\ldots <T_{\nu } <\ldots {\kern 1pt} {\kern 1pt} \, \, {\rm with}\, \, \,  T_{\nu } \to \infty  \end{equation} 
there exist a neighborhood $\Pi $ of zero with $B[0,R]:=\left\{x\in {\mathbb R}^{ n} :|x|\le R\right\}\subset \Pi $ and a map $k:{\mathbb R}^{+} \times \Pi \to  {\mathbb R}^{ m} $ such that for any $x\in \Pi $ the map $k(\cdot ,x): {\mathbb R}^{+} \to  {\mathbb R}^{ m} $ is measurable and locally essentially bounded and the trajectory $x(\cdot )$ of the sampled-data closed loop system
\begin{align}   \dot{x}=f(x,k(t,x(T_{i} &))),\, \, t\in [T_{i} ,\; T_{i+1} ),\; \, \, i=1,2,\ldots \nonumber \\  & x(0)\in \Pi   \label{2.9} \end{align} 
satisfies:
\begin{flalign}\text{Stability:} &&\hspace{-2em}
\begin{array}{c}\forall \varepsilon >0\Rightarrow \exists \delta =\delta (\varepsilon )>0:x(0)\in \Pi , \\
|x(0)|\le \delta\Rightarrow |x(t)|\leq \varepsilon ,\; \forall t\geq 0\end{array}
&\label{2.10}
\end{flalign}
\begin{flalign}\text{Attractivity:}  && \mathop{\lim }\limits_{t\to \infty } x(t)=0,\; \forall x(0)\in \Pi&&\label{2.11}\end{flalign}
\end{dfn}
\vspace{1em}

It should be pointed out that Definition 2 is stronger than the concept of sampled-data semi-global asymptotic stabilization adopted in earlier relative works in the literature, because the partition of times in \eqref{2.8} is arbitrary.
  
The proof of the following proposition is based on a  generalization of the methodology applied in \cite{art:17} and is provided in \cite{art:arX}:\\[-0.4em]

\begin{prop}\label{prp:2}
Under Assumption 1, system \eqref{1.1} is SDF-SGAS and therefore SDF-WGAS.\\[-0.4em]
\end{prop}

We next present the statement of the central result of present work, which provides a Lie algebraic sufficient condition for SDF-SGAS(WGAS) for the affine in the control single-input system \eqref{1.2}. In the sequel we assume that its dynamics $f$, $g$ are smooth ($C^{\infty}$). We denote by $ Lie\{ f,g\}$ the Lie algebra generated by $\{f,g\}$. Also, let $L_{1}:=span\{f,g\}$ and $L_{i+1}:=span\{[X,Y],\, X\in L_{i},Y\in L_{1}\}$, $i=1,2,\ldots$ and for any nonzero $\Delta \in Lie\{f,g\}$ define
\begin{equation} \label{2.12} order_{\left\{f,g\right\}}  \Delta 
\left\{\hspace{-0.2em}\begin{array}{l} 
{:=1,\, \mathrm{if}\; \Delta \in L_{1}\setminus\{0\} } \\ 
{\hspace{-0.6em}\left.\begin{array}{l}:=k>1,\mathrm{if}\, \Delta=\Delta_{1}+\Delta_{2},\,\\
\hspace{2em}\textrm{with}\;\Delta_{1}\in L_{k}\setminus\{0\}\,\,\textrm{and} \\
\hspace{2em}\Delta_{2}\in span\{\cup_{i=1}^{i=k-1} L_{i}\}
\end{array}\right.} 
\end{array} \right. 
 \end{equation} 
As a consequence of Proposition 2 we get:\\[-0.4em]
\begin{prop}\label{prp:3}
For  the affine in the control case \eqref{1.2} assume that there exists a smooth function $V: {\mathbb R}^{ n} \to {\mathbb R}^{+} $, being positive definite and proper, such that for every $x\neq 0$, either $(gV)(x)\neq 0$, or one of the following properties hold:\\ Either
\begin{equation}  (gV)(x )=0\Rightarrow (fV)(x )<0  \label{2.13}\end{equation}
or there exists an integer $N=N(x)\ge 1$ such that
\begin{subequations} \label{2.14}
\begin{equation}
\label{2.14a}
(gV)(x)=0,\; (f^{i} V)(x)=0,\; \; i=1,2,\ldots ,N\end{equation}
\begin{align}
(\Delta_{{1} }& \Delta_{{2} } \ldots \Delta_{{k} } V)(x)=0\nonumber \\
\forall \Delta _{{1} } , \Delta_{{2} }& ,\ldots , \Delta_{{k} } \in Lie\{ f,g\} \setminus \{ g\}\nonumber\\
\mathrm{with}&\, \, \sum _{p=1}^{k}order_{\{f,g\}} \Delta_{{p} }  \le N \label{2.14b}
\end{align}
\end{subequations}
where $(f^iV)(x):=f(f^{i-1}V)(x)$, $i=2,3,\ldots$, $(f^1V)(x):=(fV)(x)$ and in such a way that one of the following properties hold:
\begin{flalign}\text{(P1)}  && (f^{N+1} V)(x)<0 && \label{2.15}\end{flalign}
(P2) $N$ is odd and
\begin{equation} \label{2.16} ([[\ldots [[f,\underbrace{g],g],\ldots ,g],g]}_{N\, \, \, \, }V)(x)\neq 0 \end{equation} 
(P3) $N$ is even and 
\begin{equation} \label{2.17} ([[\ldots [[f,\underbrace{g],g],\ldots ,g],g]}_{N\, \, \, \, }V)(x)<0 \end{equation} 
(P4) $N$ is an arbitrary positive integer with
\begin{subequations}\label{2.18}
\begin{align}(f^{N+1} V)(x)&=0,\;    \label{2.18a}\\
([[\ldots [[g,\underbrace{f],f],\ldots ,f],f]}_{N\, \, \, \, }&V)(x)\ne 0\; \label{2.18b}
 \end{align}
\end{subequations}
Then system \eqref{1.2} is SDF-SGAS and therefore SDF-WGAS.
\end{prop}

\begin{rem} For the particular case of $N=1$, condition \eqref{2.14a} is equivalent to $(gV)(x)=0$ and $(fV)(x)=0$, the previous equality is equivalent to \eqref{2.14b} and obviously \eqref{2.16} is equivalent to $([f,g]V)(x)\ne 0$. It turns out, according to the statement of Proposition 3, that, under \eqref{1.3}, the system \eqref{1.2} is SDF-SGAS and therefore SDF-WGAS; the latter conclusion, namely, that (1.3) implies SDF-WGAS, is the precise statement of \cite[Proposition 2]{art:19}. 
\end{rem}

An interesting consequence of Proposition 3 concerning 3-dimensional systems \eqref{1.2} is the following result:\\
\begin{corollary}
Consider the 3-dimensional system \eqref{1.2} and assume that:
\begin{flalign}\text{(I)} && span\{g(x), [f,g](x), [f,[f,g]](x)\}=\mathbb{R}^{3}&& \label{dimbrocket}\end{flalign}
(II) There exists a smooth positive definite and proper function $V:\mathbb{R}^{n}\rightarrow\mathbb{R}^{+}$ such that \begin{equation}\label{2.20}DV(x)\neq 0, \, \, \forall x\neq 0\end{equation} and in such a way that, either \eqref{2.13} holds, or 
\begin{equation}
(gV)(x)=0\Rightarrow (f^iV)(x)=0,\,\forall x\neq 0,\, i=1,2,3\label{2.21}
\end{equation} 
Then the system is SDF-SGAS.\\[-0.4em]
\end{corollary}

We finally consider the following interesting case of 3-dimensional systems: 
\begin{align} 
\dot{x}_1=a(x_1,x_2,x_3)x_3^L,\,\,
&\dot{x}_2=b(x_1,x_2,x_3)x_3,\,\,
\dot{x}_3=u,\,\,\nonumber\\
(&x_1,x_2,x_3)\in\mathbb{R}^{3}\label{int:e:2}
\end{align}
where 
\begin{equation}\label{H1:C2}L\geq 3\,\, \mathrm{is\,\, a\,\, positive\,\,odd\,\, integer}\end{equation} and the functions $a,b:\mathbb{R}^3\rightarrow\mathbb{R}$ are smooth ($ C^{\infty}$) and  satisfy
\begin{equation}\label{H2:C2}
a(x), b(x)\neq 0,\,\, \forall x\in \mathbb{R}^3
\end{equation}
It can be easily verified that \eqref{int:e:2} does not satisfy the well known Brockett's condition for smoothly static feedback stabilization. For $a(\cdot)=b(\cdot)=1$, it was established in \cite{art:7.01}, that \eqref{int:e:2} is small time locally controllable and in \cite{art:Sep} that is \textit{locally} asymptotically stabilizable by means of a smooth time-varying periodic feedback. We use the result of Proposition 3 of present work to establish the following result.
\begin{corollary}
Under hypotheses \eqref{H1:C2} and \eqref{H2:C2}, system \eqref{int:e:2} is SDF-SGAS.
\end{corollary}
\section{Proofs}
\noindent\textit{Outline of proof of Proposition 3:} (As mentioned, the complete proof is found in \cite{art:arX})
 Let $0\ne x_{0} \in {\mathbb R}^{ n} $ and suppose first that, either $(gV)(x_{0} )\ne 0$, or \eqref{2.13} is fulfilled, namely, $(gV)(x_{0} )=0$ and $(fV)(x_{0} )<0$. Then there exists a constant input $u$ such that both \eqref{2.7a} and \eqref{2.7b} hold; particularly, for every sufficiently small $\varepsilon >0$ we have: 
\begin{equation} \label{3.17} V(x(s,0,x_{0} ,u))<V(x_{0} ),{\kern 1pt} \forall s\in (0,\varepsilon ] \end{equation} 
Assume next that there exists an integer $N=N(x_{0} )\geq 1$ satisfying \eqref{2.14}, as well as one of the properties (P1), (P2), (P3), (P4) with $x=x_0$. Then \eqref{2.14} implies: 
\begin{equation}\label{3.10:cdc}
(fV)(x_0)=(gV)(x_0)=0
\end{equation}
In order to derive the desired conclusion, we proceed as follows. Define: 
\begin{equation} \label{3.18} X:=f+u_{1} g,\, \, Y:=f+u_{2} g \end{equation} 
and let us denote by $X_{t} (z)$ and $Y_{t} (z)$ the trajectories of the systems $\dot{x}=X(x)$ and $\dot{y}=Y(y)$, respectively, initiated at time $t=0$ from some $z\in  {\mathbb R}^{n} $. Also, for any constant $\rho>0$ define:
\begin{subequations}\label{3.19A}\begin{align}  R(t):=(X_{\rho t} \circ Y_{t} )&(x_{0} ),\,t\ge 0,\,R(0)=x_{0}   \label{3.19}\\
 m(t):&=V(R(t)),t\ge 0  \label{3.20}\end{align}\end{subequations}
and denote in the sequel by $\mathop{m}\limits^{(\nu )} (\cdot ),\, \nu =1,2,...$ its $\nu $-time derivative.  By taking into account \eqref{3.10:cdc}-\eqref{3.19A} and exploiting the Campbell-Baker-Hausdorff formula for the right hand side map of \eqref{3.19}, together with an induction procedure, it can be shown that
\begin{subequations} 
\begin{equation} \label{3.26} \mathop{m}\limits^{(1)} (0)=0 \end{equation}
and for every integer $n\geq 2$, the $n$-time derivative $\mathop{m}\limits^{(n)} (\cdot)$ of $m(\cdot)$ satisfies
\begin{align}&\mathop{m}\limits^{(n)} (0) \in (A_{0}^{n} V)(x_{0} )\nonumber\\
&+span\left\{\hspace{-0.3em}\begin{array}{l} {\rho^{r_{n} } (A_{i_{1} } A_{i_{2} } ...A_{i_{\nu } } V)(x_{0} ):{\kern 1pt} \nu \ge 2;}\\{i_{1} ,i_{2} ,...i_{\nu } \in {\mathbb N}_{0} ;\sum _{j=1}^{\nu } order_{\left\{X,Y\right\}} A_{i_{j} } =n;} \\ {r_{n} =\sum _{j=1}^{\nu } i_{j} \in \left\{1,2,...,n-2\right\}} \end{array}\right\}\nonumber\\
&+\rho^{n-1} (A_{n-1} V)(x_{0} )\label{3.34} \end{align}
\end{subequations}
where
\begin{equation} \label{3.22} 
\begin{array}{l} A_{0}   :=  {\rho X+Y,} \\ A_{\nu } := {[...[[Y,\underbrace{X],X],\ldots ,X]}_{\nu },\nu =1,2,...} 
\end{array} 
\end{equation} 
Since $A_{\nu } \in Lie\{ X,Y\} $, we may define, according to \eqref{2.12}, the order of each $A_{\nu } $ with respect to the Lie algebra of $\{ X,Y\} $; particularly, in our case, we have:
\begin{equation} \label{3.23} order_{_{\left\{X,Y\right\}} } A_{\nu } =\nu +1,\, \,  \forall \nu =0,1,2,\ldots  \end{equation} 
By taking into account definition \eqref{3.18} of the vector fields $X$ and $Y$ and by setting
\begin{equation} \label{3.35} u_{2} =-\rho u_{1}, \, \rho>0  \end{equation} 
we get
\begin{align} 
&A_{0} =(\rho+1)f, \, \, A_{1} =(\rho+1)u_{1} [f,g],\nonumber\\
&A_{2} =  (\rho+1)(u_{1}^{2} [[f,g],g]-u_{1} [[g,f],f])\nonumber \\
&\hspace{2em}\vdots \nonumber\\
&A_{n} = (\rho+1)u_{1}^{n} [... [[f,\underbrace{g],g],... ,g]}_{n{\kern 1pt} {\kern 1pt}}\nonumber\\
&+(\rho+1)u_{1}^{n-1} ([[[... [f,\underbrace{g],... ,g],g}_{n-1{\kern 1pt} {\kern 1pt} }],f]\nonumber \\
  & +[[[... [f,\underbrace{g],... ,g}_{n-2{\kern 1pt} {\kern 1pt} }],f],g]+... \nonumber\\
  &+[... [[[f,g],f],\underbrace{g]... ,g]}_{n-2{\kern 1pt} {\kern 1pt} })+[... [[[f,g],f],\underbrace{g]... ,g]}_{n-2{\kern 1pt} {\kern 1pt} })\nonumber\\
&+... +(\rho+1)u_{1}^{2} ([[[... [[f,g],\underbrace{f],... ,f],f}_{n-2{\kern 1pt} {\kern 1pt} }],g]\nonumber\\
 &+[[[... [[f,g],\underbrace{f],... ,f}_{n-3{\kern 1pt} {\kern 1pt} }],g],f]\nonumber\\
  &+... +[[... [[[f,g],g],\underbrace{f]... ,f],f}_{n-2{\kern 1pt} {\kern 1pt} }])\nonumber \\
 &-(\rho+1)u_{1} [...[[g,\underbrace{f],f],... ,f}_{n{\kern 1pt} {\kern 1pt} }],\; n=3,4,...  \label{3.36}
\end{align}
Obviously, \eqref{3.36} implies: 
\begin{align}
A_{k} \in span\{\Delta \in Lie\left\{f,g\right\}\setminus \{g\}&:order_{\{f,g\}}\Delta  =k+1\} \nonumber\\
&k=0,1,2,\ldots \label{3.37} 
 \end{align} 
Also, we recall from \eqref{3.34} and \eqref{3.23}  that 
$r_{n} =\sum _{s=1}^{\nu } i_{s} \in \{ 1,2,\ldots ,n-2\} $ and $\sum _{j=1}^{\nu } order_{\{X,Y\}} A_{i_{j}}=r_{n}+\nu =n$ with $\nu \ge 2$
and  therefore $\nu \le n-1$. By \eqref{3.34}-\eqref{3.37} and the previous facts we get:
\begin{align} \mathop{m}\limits^{(n)} &(0)\in (\rho+1)^{n} (f^{n} V)(x_{0} )+u_{1}^{} \pi _{1} (\rho,\rho+1;x_{0} )\nonumber \\
 &+span\left\{u_{1}^{k} \pi _{k} (\rho,\rho+1;x_{0} ),k=2,...,n-2\right\}\nonumber\\
&+\rho^{n-1} (\rho+1)u_{1}^{n-1} ([\ldots [[f,\underbrace{g],g],\ldots ,g}_{n-1\, \,  }]V)(x_{0} ) \nonumber \\ 
&-\rho^{n-1} (\rho+1)u_{1} ([\ldots [[g,\underbrace{f],f],\ldots ,f}_{n{\kern 1pt} -1\, \, \, \, {\kern 1pt} }]V)(x_{0} ) \label{3.38} 
\end{align} 
for $n=2,3,...$ and for certain smooth functions $\pi _{k} : {\mathbb R}^{2} \times {\mathbb R}^{n} \to  {\mathbb R},\,\,  k=1,2,\ldots ,n-2$ satisfying the following properties:

\noindent(S1) For every $x_{0} \in  {\mathbb R}^{n} $, each map $\pi _{k} (\alpha ,\beta ;x_{0} ): {\mathbb R}^{2} \to  {\mathbb R}$ is a polynomial with respect to the first two variables  in such a way that
\begin{equation} \label{3.39} \begin{array}{l} { span\{ \pi _{k} (\alpha ,\beta ;x_{0} ),\, {\kern 1pt} k=1,2,\ldots ,n-2\} \subset } \\ {span\{ (\Delta _{{1} }  \Delta _{{2} } ...\Delta _{{i} } V)(x_{0} );\; {i} \in {\mathbb N},\; } \\ {\quad {\kern 1pt} \; \; \; \Delta _{{1} }  ,\Delta _{{2} } ,...,\Delta _{{i} } \in Lie\{ f,g\} \backslash \{ g\} ;}\\{ \quad \quad \quad \sum _{j=1}^{j=i} order_{\left\{f,g\right\}} \Delta _{{j} } =n\, \, \} } \end{array} \end{equation} 
(S2) For each $x_{0} \in  {\mathbb R}^{n} $ there exist integers $\lambda _{i}$, $\mu _{i}$, $i=1,2,...,\; L\in {\mathbb N}$ with $1\le \lambda _{i} \le n-2$, $2\le \mu _{i} \le n-1$ such that the map $\pi _{1} (\alpha ,\beta ;x_{0} ): {\mathbb R}^{2} \to \mathbb{R}$ satisfies:
$$ \pi _{1} (\alpha ,\beta ;x_{0} )\in span\left\{\alpha^{\lambda _{1} } \beta ^{\mu _{1} } ,\alpha^{\lambda _{2} } \beta ^{\mu _{2} } ,...,\alpha^{\lambda _{L} } \beta ^{\mu _{L} } \right\}$$ The latter implies that for each fixed $x_{0} \in  {\mathbb R}^{n} $ the polynomials $\pi _{1} (\rho,\rho+1;x_{0} )$ and $$-\rho^{n-1} (\rho+1)([\ldots [[g,\underbrace{f],f],\ldots ,f}_{n -1\, \, \, \,  }]V)(x_{0} )$$ are linearly independent, provided that 
\begin{equation} \label{3.41} ([[\ldots [[g,\underbrace{f],f],\ldots ,f],f]}_{n-1\, \, \, \, }V)(x_{0} )\ne 0\;  \end{equation} 
If we define:
\begin{align} \label{3.42} \xi _{n} (\rho;x):=&\pi _{1} (\rho ,\rho+1;x_{0} )\\
&-\rho^{n-1} (\rho+1)([\ldots [[g,\underbrace{f],f],\ldots ,f}_{n{\kern 1pt} -1\, \, \, \, {\kern 1pt} }]V)(x_{0} )\nonumber \end{align} 
the inclusion \eqref{3.38} is rewritten:    
\begin{align} &\mathop{m}\limits^{(n)} (0)\in (\rho+1)^{n} (f^{n} V)(x_{0} )+u_{1}^{} \xi _{n} (\rho;x_{0} )\nonumber \\
&+span\left\{u_{1}^{k} \pi _{k} (\rho,\rho+1;x_{0} ),k=2,...,n-2\right\}\nonumber \\ 
&+\rho^{n-1} (\rho+1)u_{1}^{n-1} ([\ldots [[f,\underbrace{g],g],\ldots ,g}_{n-1\, \, {\kern 1pt} {\kern 1pt} }]V)(x_{0} ) \label{3.43} \end{align} 
and a  constant $\rho=\rho(x_{0} )>0$ can be found with 
\begin{equation} \label{3.44} \begin{array}{l} {\xi _{n} (\rho;x_{0} )\ne 0} \end{array} \end{equation} 
provided that \eqref{3.41} holds. Suppose now that there exists an integer $N=N(x_{0} )\ge 1$ satisfying  \eqref{2.14}, as well as one of the properties  (P1), (P2), (P3), (P4) with $x=x_0$. By \eqref{3.26} and by taking into account \eqref{2.14}, \eqref{3.38} and \eqref{3.39} it  follows:
\begin{equation} \label{3.45} \mathop{m}\limits^{(n)} (0)=0,{\kern 1pt} {\kern 1pt} n=1,2,\ldots ,N \end{equation} 
and also, by taking into account \eqref{3.19}, \eqref{3.20} and \eqref{3.35}-\eqref{3.44} it can be shown that in all cases considered in the statement of Proposition 3, there exist constants $\rho=\rho(x_0)>0$ and $u_1$ such that, if we define
\begin{equation} \label{3.47} w(s;t,x_{0}):=\left\{\begin{array}{c} {u_{2} =-\rho u_{1} ,{\kern 1pt} s\in [0,t]} \\ {u_{1} ,{\kern 1pt} s\in (t,t+\rho t]} \end{array}\right.  \end{equation} 
it holds $\mathop{m}\limits^{(N+1)}(0)<0$, which, in conjunction with \eqref{3.45}, asserts that for every sufficiently small $\sigma=\sigma(x_0) >0$ we have 
\begin{subequations}\begin{equation}m(t)<m(0),\,\, \forall t\in (0,\sigma ]\label{mt:L:m0}\end{equation} where 
\begin{align}m(t):&=V((X_{\rho t} \circ Y_{t} )(x_{0} ))\nonumber\\
&=V(x(t+\rho t,0,x_{0} ,w(\cdot;t,x_{0} ))\label{mV}\end{align} \end{subequations}
and $x(\cdot ,0,x_{0} ,w(\cdot;t,x_{0} ) )$ is the trajectory of \eqref{1.2} corresponding to the input $w(\cdot;t,x_{0} ) $. Equivalently:
\begin{equation} \label{3.48} V(x(t,0,x_{0} ,w(\cdot;\tfrac{t}{1+\rho},x_{0} ) ))<V(x_{0} ){\kern 1pt} ,\forall t\in(0,\tfrac{\sigma}{1+\rho}] \end{equation} 
hence, we may pick $\varepsilon \in (0,\tfrac{\sigma}{1+\rho} ]$ sufficiently small in such a way that inequality in \eqref{3.48} holds for $t:=\varepsilon $, namely,
\begin{subequations}\label{V:p3}\begin{equation}V(x(\varepsilon ,0,x_{0} ,u(\cdot,x_{0}) ))<V(x_{0} )\end{equation} with 
$u(s,x_0):=w(s;\tfrac{\varepsilon}{1+\rho},x_0),\,\, s\in(0,\varepsilon]$
and simultaneously 
\begin{equation}V(x(s,0,x_{0} ,u(\cdot,x_{0}) ))\le 2V(x_{0} ),{\kern 1pt} {\kern 1pt} \forall s\in (0,\varepsilon ] \end{equation}\end{subequations}
We conclude, by taking into account \eqref{3.17} and \eqref{V:p3}, that for every $x_{0} \ne 0$ and $\xi>0$, there exist $\varepsilon=\varepsilon(x_0)\in(0,\xi]$ and a measurable and locally essentially bounded control $u(\cdot,x_{0}) :[0,\varepsilon ]\to {\mathbb R}$ such that \eqref{2.7a} and \eqref{2.7b} hold with $a(s):=2s$. Therefore, according to Proposition 2, \eqref{1.2} is SDF-SGAS.    $\blacksquare$

\noindent\textit{Proof of Corollary 1:} 
First, by invoking assumptions \eqref{dimbrocket} and \eqref{2.20} it follows that for every $x\neq 0$, either $(gV)(x)\neq 0$, or
\begin{equation}\label{gv:eq:0}
(gV)(x)=0
\end{equation}
 which in conjunction with \eqref{2.13} implies the desired statement. Also, by virtue of \eqref{dimbrocket}-\eqref{2.21}, we have
\begin{subequations}
\begin{align}
(fV)(x)=(f^2V)(x)=(f^3V)(x)=0\label{fV:eq:0}\\
|([f,g]V)(x)|+|([f,[f,g]]V)(x)|\neq 0\label{sum:brackets}
\end{align}
\end{subequations}
For those $x\neq 0$ for which \eqref{gv:eq:0} holds, we consider two cases. The first is $([f,g]V)(x)\neq 0$, which in conjunction with \eqref{gv:eq:0} and \eqref{fV:eq:0} assert that \eqref{2.14a} and (P4) hold with $N=1$. The other case is
\begin{subequations}\label{bracket:2}
\begin{align}
([f,g]V)(x)=0\label{fgV:eq:0}\\
([f,[f,g]]V)(x)\neq 0 \label{ffgV:neq:0}
\end{align}
\end{subequations}
which in conjunction with \eqref{gv:eq:0} and \eqref{fV:eq:0} assert that \eqref{2.14a}, \eqref{2.14b} and (P4) are fulfilled with $N=2$. We conclude, according to the statement of Proposition 3, that the 3-dimensional system \eqref{1.2} is SDF-SGAS.$\blacksquare$\\
 
\noindent\textit{Proof of Corollary 2:} 
We define:
\begin{equation}\begin{array}{c}f(x):=(a(x_1,x_2,x_3)x_3^L,\,b(x_1,x_2,x_3)x_3,0)^T,\\
g(x):=(0,0,1)^T,\,\,x:=(x_1,x_2,x_3)^T\end{array}\label{e:dyn}
\end{equation}
and 
\begin{equation} 
V(x):=\tfrac{1}{2}x_1^2+\tfrac{1}{L+1}x_2^{L+1}+\tfrac{1}{2}x_3^2\label{lyap}\end{equation}
that obviously is positive definite and proper. According to the previous definitions, it follows that 
\begin{subequations}\label{aa}
\begin{equation}
([f,g])(x)=\left(\begin{matrix}-\frac{\partial a}{\partial x_3}(x_1,x_2,x_3)x_3^L-La(x_1,x_2,x_3)x_3^{L-1}\\-\frac{\partial b}{\partial x_3}(x_1,x_2,x_3)x_3-b(x_1,x_2,x_3)\\0\end{matrix}\right)\label{a1}
\end{equation}
and for each integer $k: 2\leq k\leq L$ it holds:
\begin{align}
&([...[[f,\underbrace{g],g]...g}_{k\,\,}])(x)=(A_{1,k}(x_1,x_2,x_3)\nonumber\\
&\hspace{2em}+(-1)^k\prod_{i=0}^{k-1}(L-i)a(x_1,x_2,x_3)x_3^{L-k},\nonumber\\
&\hspace{2em}A_{2,k}(x_1,x_2,x_3)+(-1)^kk\frac{\partial^{k-1} b}{\partial x_3^{k-1}}(x_1,x_2,x_3),0)^T\label{a2}\\
&([...[[g,\underbrace{f],f]...f}_{k\,\,}])(x)=(B_{1,k}(x_1,x_2,x_3),B_{2,k}(x_1,x_2,x_3),0)^T\label{a3}
\end{align}\end{subequations}
for certain smooth functions $A_{1,k},A_{2,k},B_{1,k},B_{2,k}:\mathbb{R}^3\rightarrow \mathbb{R}$, satisfying 
\begin{subequations} \label{hqps:all}
\begin{equation}
A_{1,k}(\cdot,\cdot,0)=A_{2,k}(\cdot,\cdot,0)=B_{1,k}(\cdot,\cdot,0)=B_{2,k}(\cdot,\cdot,0)=0,\label{hqps}
\end{equation}
and
\begin{align} 
&\frac{\partial^{j} A_{1,n}}{\partial x_2^j}(\cdot,\cdot,0)=\frac{\partial^j B_{1,n}}{\partial x_2^j}(\cdot,\cdot,0)=\frac{\partial^j B_{2,n}}{\partial x_2^j}(\cdot,\cdot,0)=0,\nonumber\\
&j=1,...,L-1;\,\,n=2,...,L-j+1;\,\, \label{hqps:der}
\end{align}
\end{subequations}
From \eqref{e:dyn}-\eqref{aa} we also get
\begin{subequations} 
\begin{align} 
&(gV)(x)=x_3\label{gV};\\
&([f,g]V)(x)=-\tfrac{\partial a}{\partial x_3}(x_1,x_2,x_3)x_1x_3^L\nonumber\\
&\hspace{6em}-La(x_1,x_2,x_3)x_1x_3^{L-1}\nonumber\\
&\hspace{2em}-\tfrac{\partial b}{\partial x_3}(x_1,x_2,x_3)x_2^Lx_3-b(x_1,x_2,x_3)x_2^L,\,\forall x\in \mathbb{R}^{3}\label{brk:fg:V1}\end{align}
and for any integer $k$: $2\leq k\leq L$ it holds:
\begin{align}
&([...[[f,\underbrace{g],g]...g}_{k\,\,}]V)(x)=A_{1,k}(x_1,x_2,x_3)x_1\nonumber\\
&\hspace{1em}+(-1)^k\prod_{i=0}^{k-1}(L-i)a(x_1,x_2,x_3)x_3^{L-k}x_1\nonumber\\
&\hspace{1em}+A_{2,k}(x_1,x_2,x_3)x_2^L+(-1)^kk\frac{\partial^{k-1} b}{\partial x_3^{k-1}}(x_1,x_2,x_3)x_2^L\label{a}\\
&([...[[g,\underbrace{f],f]...f}_{k\,\,}]V)(x)=B_{1,k}(x_1,x_2,x_3)x_1\nonumber\\
&\hspace{10em}+B_{2,k}(x_1,x_2,x_3)x_2^L\label{b}
\end{align}
\end{subequations}
Let $x\neq0$ for which
\begin{equation} (gV)(x)=x_3=0\label{gv:0}\end{equation}
It then follows by virtue of \eqref{e:dyn}, \eqref{lyap} and \eqref{gV} that \begin{equation}\label{fV} (f^kV)(x)=0,\,\,k=1,2,...\end{equation} therefore (2.14a) holds, and further, by invoking \eqref{H1:C2}
\begin{equation}
([f,g]V)(x)=-b(x_1,x_2,0)x_2^L\label{brk:fg:V}
\end{equation} 
Then we may distinguish the following two cases:

\noindent\textbf{Case 1:} $x_2\neq0$ with $x_1\neq0$ and $x_3=0$. Then by taking into account our hypothesis (2.24) and \eqref{brk:fg:V}, it follows that $([f,g]V)(x)\neq 0$, which in conjunction with \eqref{fV} asserts that both (2.14) and (P2) in the statement of Proposition 3 are satisfied with $N=1$.

\noindent \textbf{Case 2:} $x_2=0$ with $x_1\neq 0$  and $x_3=0$. It then follows from \eqref{H2:C2}, \eqref{a}, \eqref{gv:0} and \eqref{brk:fg:V} that  
\begin{subequations}\label{brk:fgg:gff:k:v}\begin{align}
&([[[f,\underbrace{g],g],...,g}_{k\,\, }]V)(x)=0, \,\,k=1,...,L-1;  \label{C2:a}
\end{align}
\begin{equation}
([...[[f,\underbrace{g],g],...,g}_{L \, \, }]V)(x)\neq0, \forall x_1\neq 0\label{C2:c}
\end{equation}\end{subequations}
and therefore we can easily verify from our hypotheses \eqref{H1:C2}, \eqref{H2:C2} and \eqref{C2:c}, that (P2) holds with $N=L$. By taking into account \eqref{fV}, it also follows that (2.14a) holds for $k=1,...,L$, thus, in order to verify that all statements of Proposition 3 are satisfied, it remains to show that (2.14b) holds as well. Particularly, we show that, if we define      
\begin{align}
&\pi _{k} (x):=(\Delta _{1} \Delta _{2} ...\Delta _{k} V)(x);   \nonumber\\
\Delta _{1} ,...,&\Delta _{k} \in Lie\{ f,g\} \backslash \{ g\}\,\, \mathrm{with}\,\, \sum _{p=1}^{k}order_{\{f,g\}}\Delta _{p} \le L\label{delta:ell}
\end{align}
it holds 
\begin{equation} \pi _{k} (x_{1} ,0,0)=0,{\kern 1pt} {\kern 1pt} \forall x_{1} \in {\rm {\mathbb R}} \label{pol:ell:0} \end{equation} 

In order to establish \eqref{pol:ell:0}, it suffices to consider in \eqref{delta:ell} only those $\Delta_p$ satisfying
$$\Delta_p\in\{f,\,\,[...[[f,\underbrace{g],g],...,g}_{k_1 \, \, }],\,\,[[g,\underbrace{f],f],...,f}_{k_2\,\, }]\}$$ for certain appropriate $k_1,k_2\in \mathbb{N}$. Notice first that, due to \eqref{e:dyn}, \eqref{aa} and \eqref{hqps:all}, each $\Delta _{p} $, $p=1,2,...,k$ above is written as 
\begin{subequations}\begin{equation} \label{Delta} \Delta _{p}(x) =(C_{1,k} (x_{1},x_2 ,x_{3} ),C_{2,k} (x_{1},x_2 ,x_{3} ),0)^{T}  \end{equation} 
for all $x=(x_1,x_2,x_3)\in\mathbb{R}^{3}$ and for certain smooth functions $C_{1,k} (\cdot ,\cdot, \cdot )$ and $C_{2,k} (\cdot ,\cdot, \cdot )$ with 
\begin{align}
&C_{1,k} (\cdot,\cdot,0)=0;\;\;\frac{\partial^{j}C_{1,q}}{\partial x^j_2}(\cdot,\cdot,0)=0,\nonumber\\
&\,\,\,j=1,...,k;\,\,q=1,...,k-j+1,\nonumber\\
&\,\,\,\mathrm{for\,\, the\,\, case}\nonumber\\
&\,\,\,\,\,\Delta _{p} \in D:=\{[...[[f,\underbrace{g],g],...,g}_{n \, \, }],n=1,...,L\}\label{D:fg:B}\end{align}and
\begin{align}
&C_{1,k} (\cdot,\cdot ,0)=0;\,\, C_{2,k} (\cdot,\cdot ,0)=0;\;\;\frac{\partial^{j}C_{1,q}}{\partial x^j_2}(\cdot,\cdot,0)=0,\nonumber\\
&\,\,\,\,\,j=1,...,k;\,\,q=1,...,k-j+1,\nonumber\\
&\,\,\mathrm{for\,\, those}\,\,\Delta _{p} \in Lie\{ f,g\} \backslash \{ \{g\}\cup D\} \label{D:B}
\end{align}
\end{subequations}
We then may use the previous facts, together with \eqref{e:dyn}-\eqref{hqps:all} and an elementary induction procedure, in order to establish that for every integer $k\in \{1,...,L-1\}$, for which the inequality in \eqref{delta:ell} holds, there exist smooth functions $\Xi_{1} =\Xi_{1} (x_{1} ,x_{2} ,x_{3} )$ and $\Xi_{2} =\Xi_{2} (x_{1},x_2 ,x_{3} )$ in such a way that  
\begin{subequations}
\begin{align}
\Xi_{1} (\cdot ,\cdot ,0)&=0\label{B1}\\
\pi _{k} (x_{1} ,x_{2} ,x_{3} )&=\Xi_{1} (x_{1},x_2 ,x_{3} )+\Xi_{2} (x_{1},x_2 ,x_{3} )x_{2}^{L-k+1}\label{pi:ell}
\end{align}
\end{subequations}
and the latter establishes \eqref{pol:ell:0}. It follows from (2.24), \eqref{fV}, \eqref{C2:c} and \eqref{pol:ell:0}  that for the Case 2, both (2.14) and (P2) hold with $N=L$.   
 
We conclude, that in both Cases 1 and 2, hypothesis of Proposition 3 is satisfied, therefore system  is SDF-SGAS.$\blacksquare$

\end{document}